\documentclass{amsart}
\usepackage{amssymb}
\usepackage{enumitem}
\usepackage{xcolor}
\usepackage{pinlabel}
\usepackage{hyperref}
\usepackage{thm-restate}

\newtheorem {theorem}{Theorem}[section]
\newtheorem {lemma} [theorem] {Lemma}
\newtheorem {proposition} [theorem] {Proposition}
\newtheorem {corollary} [theorem] {Corollary}

\newtheorem* {claim*}{Claim}

\newtheorem* {subclaim*}{Subclaim}

\theoremstyle{definition}
\newtheorem{remark}[theorem]{Remark}
\newtheorem {definition} [theorem] {Definition}
\newtheorem {example} [theorem] {Example}

\newcommand{\mc}{\mathcal}
\newcommand{\from}{\colon\thinspace}

\newcommand{\Homeo}{\operatorname{Homeo}}

\newcommand{\Hom}{\operatorname{Hom}}

\newcommand{\dvis}{d_\partial}
\newcommand{\dgam}{d_{X}}

\newcommand{\bgam}{\partial \Gamma}
\newcommand{\diam}{\operatorname{diam}}

\newcommand{\idgamma}{e}

\newcommand{\minus}{-}
\newcommand{\bN}{\mathbb{N}}
\newcommand{\bH}{\mathbb{H}}
\newcommand{\bR}{\mathbb{R}}
\newcommand{\bZ}{\mathbb{Z}}
\newcommand{\trot}{\operatorname{\tilde{r}ot}}

\newcounter{notes}

\newcommand{\putinbox}[1]{\noindent\fbox{\parbox{.98\textwidth}{#1}}}

\title{Stability of hyperbolic groups acting on their boundaries}
\author[K. Mann]{Kathryn Mann}
\address{Department of Mathematics, Cornell University, Ithaca, NY 14853, USA
}
\email{k.mann@cornell.edu}

\author[J.F. Manning]{Jason Fox Manning}
\address{Department of Mathematics, Cornell University, Ithaca, NY 14853, USA
}
\email{jfmanning@cornell.edu}

\author[T. Weisman]{Theodore Weisman}
\address{Department of Mathematics, University of Michigan, Ann Arbor,
  MI 48109, USA
}
\email{tjwei@umich.edu}

\begin{document}

\begin{abstract}
  A hyperbolic group $\Gamma$ acts by homeomorphisms on its Gromov boundary $\bgam$.  We use a dynamical coding of boundary points to show that such actions are {\em topologically stable} in the dynamical sense: any nearby action is semi-conjugate to (and an extension of) the standard boundary action.

  This result was previously known in the special case that $\bgam$ is a topological sphere.  Our proof here is independent and gives additional information about the semiconjugacy in that case.  Our techniques also give a new proof of global stability when $\bgam = S^1$.
\end{abstract}

\maketitle

\section{Introduction} 
A discrete, hyperbolic group $\Gamma$, viewed as a (coarse) metric space, admits a natural compactification by its
{\em Gromov boundary}, denoted $\bgam$.  This boundary is a compact
metrizable space, and if $\Gamma$ is not virtually cyclic, it has
no isolated points. The action of $\Gamma$ on itself by left-multiplication extends naturally to 
an action on $\bgam$ by homeomorphisms, with rich global dynamics.  In 
fact, Bowditch \cite{Bowditch98} showed that a certain dynamical property (convergence group dynamics with all conical limit points) {\em characterizes} boundary actions of hyperbolic groups among all actions of all groups on perfect, compact, metrizable spaces.

This paper considers a different basic dynamical question with a long history, namely stability under perturbation.  
The general study of stability or rigidity of actions of groups on boundaries dates back to Mostow and Furstenburg, who worked in an algebraic rather than dynamical context.   On the dynamical side, Sullivan \cite{Sullivan} showed that the actions of convex cocompact Kleinian groups on their limit sets are stable under $C^1$--small perturbation; his techniques were generalized to prove stability under ``Lipschitz-small" perturbations for a much broader class of group actions, including those of hyperbolic groups on their boundaries, in \cite{KKL}.

Here we consider the more general question of $C^0$--small perturbations and stability in the sense of topological dynamics.  
In this setting, a representation $\rho_0 \in \Hom(\Gamma, \Homeo(X))$ is said to be a {\em topological factor} of another such representation $\rho$ if there is a continuous, surjective map $h: X \to X$ such that $h \circ \rho(\gamma) = \rho_0(\gamma) \circ h$ for all $\gamma \in \Gamma$.  In this case, $\rho$ is said to be an {\em extension} of $\rho_0$, and  $h$ is called a {\em semi-conjugacy}.  
An action $\rho$ is said to be {\em topologically} or {\em $C^0$--stable} if
every nearby action is an extension of it, or in other words, nearby actions encode the same dynamical information as $\rho$.  
We prove the following.  

\begin{theorem}[Boundary actions are $C^0$ stable] \label{thm:main} 
Let $\Gamma$ be a hyperbolic group.  
For any neighborhood $\mathcal{U}$ of the identity in the space of continuous self-maps of $\bgam$, there exists a neighborhood $\mathcal{V}$ of the standard boundary action in $\Hom(\Gamma, \Homeo(\bgam))$ such that any $\rho \in \mathcal{V}$ is an extension of the standard boundary action, via a semi-conjugacy in $\mathcal{U}$.  
\end{theorem} 
In particular, Theorem~\ref{thm:main} implies that actions
$C^0$--close to the standard boundary action cannot have larger kernel
than the standard action.  This kernel is finite as long as $\Gamma$
is not virtually cyclic.

After this paper was completed, we discovered that Gromov had
outlined a strategy for proving Theorem~\ref{thm:main} in
\cite[8.5.Y]{GromovEssay} using what he calls \emph{non-expanding
  quasi-geodesic fields}.  The reader may find it interesting to compare Gromov's suggestions with our proof.

$C^0$--stability of boundary actions was previously proved for the
special case of fundamental groups of compact Riemannian manifolds of
negative curvature in \cite{BowdenMann}, and for hyperbolic groups
with sphere boundary in \cite{MannManning}.  
Our argument here gives a stand-alone proof of the previously mentioned results, and provides additional information about the semi-conjugacy.  
\begin{corollary}\label{cor:cell}
  In case $\bgam$ is a topological $n$--sphere, one may choose a neighborhood $\mathcal{V}$ of the standard action so that all semi-conjugacies as in the conclusion of Theorem~\ref{thm:main} are cellular maps.  In particular, the semi-conjugacies can be approximated by homeomorphisms.
\end{corollary}

Our methods, combined with Ghys' work on semi-conjugacy of circle maps, also give a new proof of the following global rigidity result for groups with $S^1$ boundary, proved in \cite{Mann} and in \cite{MatsumotoBP}.  

\begin{restatable}[Global stability for $S^1$ boundary]{theorem}{globalstab}\label{thm:global} 
Let $\Gamma$ be a hyperbolic group with $\bgam = S^1$.  Let $\rho_0\from \Gamma \to \Homeo(S^1)$ be either the standard boundary action, or a lift of this action to a $k$-fold cover of $S^1$, for some $k \in \mathbb{Z}$.  Then the connected component of $\rho_0$ in $\Hom(\Gamma, \Homeo(S^1))$ is equal to the set of actions that are semi-conjugate to $\rho_0$ by a monotone, degree one semi-conjugacy. 
\end{restatable} 

The strategy of proof in both~\cite{BowdenMann} and~\cite{MannManning} was
to translate perturbations of actions into nice maps between foliated
spaces, translating the dynamical problem into a geometric one.  Here,
we use a different approach, inspired by the proof of a stability
property for relative Anosov representations given in
\cite{Weisman}.  This approach more closely follows the idea of
\emph{dynamical coding} of boundary points originally employed by
Sullivan, and involves the construction of an automaton that outputs
quasi-geodesic strings.  While the coding we construct is tailored to
proving stability of the action, the use of automata in the study of
hyperbolic groups has a long history dating back to Cannon
\cite{Cannon}; see \cite[Ch.3]{word_processing} for an introduction.

\subsection{Necessity of semi-conjugacy}
One cannot hope to improve the semi-con\-ju\-gacy in Theorem \ref{thm:main} to a genuine conjugacy without stronger restrictions on the perturbation.  One can easily build examples of arbitrarily small perturbations of the action of a free group on its boundary 
which are not conjugate to the original action, as in the following example.  

\begin{example}[Perturbation of action of $F_2$]  Let $F_2$ be the free group on the letters $a$ and $b$.  The boundary $\partial F_2$ is a Cantor set.
  We start by modifying the action of $a$ in a small clopen neighborhood $N$ of its attracting fixed point $x_+$, as follows. Pick a fundamental domain $D$ for the standard action, so that the sets $a^k(D)$, $k = 0, 1, 2 \ldots$ partition $N - \{x_+\}$ into countably many disjoint, clopen sets.  Let $X \subset N$ be a proper clopen neighborhood of $x_+$, and let $x'\in N-X$.  Let $\{D_k\}_{k=0}^\infty$ be a partition of $N-(X\cup\{x'\})$ into countably many clopen subsets of decreasing diameter accumulating to $x'$, and set $D_{-1} = a^{-1}D$.  Define a modified action of $a$ on $N$ as follows.  On $\partial F_2 - (N\cup a^{-1}D)$, the action is unchanged.  For each $k\ge -1$ define the restriction of $a$ to $D_k$ to be a homeomorphism to $D_{k+1}$.  On $X$, define $a$ to be the identity.  This gives a homeomorphism of the Cantor set which is not conjugate to the original action of $a$, since it has uncountable fixed point set.
  One may extend this action of $a$ to an action of $F_2$ by inserting the standard action of $b$ (or a similar modification if desired).  
\end{example} 

Examples of non-conjugate perturbations of actions of Kleinian groups
with sphere boundary are given in \cite{BowdenMann}.

\subsection{Relatively hyperbolic groups}
It should be possible to use our methods to prove a relative version
of Theorem~\ref{thm:main}, where $\Gamma$ is replaced by a relatively
hyperbolic group and $\bgam$ is replaced by the Bowditch boundary.  In
fact much of \cite{Weisman} (adapted in this paper) takes place in
this setting.

However, the action of a relatively hyperbolic group on its Bowditch
boundary is not stable in general.  That is, if we
replace ``hyperbolic'' with ``relatively hyperbolic'' and ``boundary''
with ``Bowditch boundary'' in Theorem~\ref{thm:main}, the statement is
no longer true.  In fact in cases where the boundary has a $C^1$
structure, the statement can even fail for $C^1$--perturbations.
As an example, consider the action of the fundamental group of a
finite volume non-compact hyperbolic 3-manifold $M$ on the ideal
boundary of $\bH^3$ (which is equivariantly identified with the
Bowditch boundary of $\pi_1M$).
By Thurston's hyperbolic Dehn filling
theorem, there are arbitrarily small $C^1$--deformations of the action
of $\pi_1M$ on $\partial \bH^3$ which have infinite kernel---and
therefore cannot be semi-conjugate to the original action, where the
kernel is trivial.
Thus any version of Theorem~\ref{thm:main} in the
relative case must in some way restrict the allowable deformations in
$\Hom(\Gamma, \Homeo(\bgam))$.

\subsection*{Outline}
Section~\ref{sec:setup} describes the process for coding boundary points.  
In Section~\ref{sec:codingproperties} we show that coding sequences give quasi-geodesics, establish a technical result describing the relationship between two codings of the same point, and discuss how conjugacy changes codings.  We also make some remarks on the relationship between our work and Bowditch's {\em annulus systems}.  
Section~\ref{sec:proof} uses the results of the previous sections to prove Theorem~\ref{thm:main} and Corollary~\ref{cor:cell}.  Finally in Section~\ref{sec:circle} we specialize to $\bgam = S^1$ and outline the proof of Theorem \ref{thm:global}

\subsection*{Acknowledgments}
K.M. was partially supported by NSF grant DMS 1844516 and a Sloan
fellowship.  J.M. was partially supported by Simons Collaboration
Grant 524176.  T.W. was partially supported by NSF grant DMS 1937215.

\section{Set-up: coding boundary points} \label{sec:setup}
We assume familiarity with the basics of hyperbolic groups and their boundaries. The reader may consult \cite[8.2]{GromovEssay} or \cite[III.H]{BH99} for a general reference.  
Let $\Gamma$ be a hyperbolic group, with fixed finite, symmetric, generating set $\mathcal{S}$.  
Theorem \ref{thm:main} is trivially true for two ended (virtually $\mathbb{Z}$) groups, so for the remainder of the paper we assume that $\Gamma$ is not virtually cyclic.  We denote by $\bgam$ the Gromov boundary of $\Gamma$, or equivalently the Gromov boundary of the Cayley graph of $\Gamma$ with respect to $\mathcal{S}$.  

Fix any metric (for instance, a visual metric) $\dvis$ on $\bgam$.
The open $\epsilon$--ball around a point $p \in \bgam$ with respect to
this metric is denoted by $B_\epsilon(p)$, and the open
$\epsilon$--neighborhood of a set $K \subset \bgam$ is denoted by
$N_\epsilon(K)$.

As is well known, the action of $\Gamma$ on the set of pairs of
distinct points in $\bgam$ is cocompact, hence we have the following.
\begin{lemma}\label{lem:distantpairs}
  There is some $D>0$ so that for every pair $a,b$ of distinct points in $\bgam$, there is a $g\in \Gamma$ so that $\dvis(ga,gb)\ge D$.
\end{lemma}
\putinbox{We fix such a constant $D$ for the rest of the paper.}

\medskip

Our first goal is to ``code'' boundary points, i.e. to associate 
to each boundary point a collection of infinite paths in a certain automaton.  In Section~\ref{sec:codingproperties} we will see that each of these paths tracks a geodesic ray limiting to the boundary point.
We wish to define the coding in a way
that uses only the dynamics of the action of $\Gamma$ on its boundary,
so that the coding sequences will still contain meaningful information
after the action of $\Gamma$ is perturbed.

The inspiration for this comes from Sullivan \cite{Sullivan}, who
codes points using sequences of elements that contract subsets of the
boundary a uniform amount.  Sullivan's ``uniform contraction'' is not
$C^0$--stable, so we instead follow the modification of this approach
given in \cite{Weisman}, and use the topological dynamics of the
action to construct nested sequences of sets that capture the idea of
uniform contraction.

We recall the following.
\begin{definition}
 Let $G$ act by homeomorphisms on a metric space $X$.  A point
$x \in X$ is a \emph{conical limit point} for the action if there exists a sequence
$\{g_n\}_{n \in \bN}$ of elements of $G$ and distinct points
$a, b \in X$ such that $g_n x \to a$ and $g_n y \to b$ for all
$y \ne x$, uniformly on compact sets in $X - \{x\}$.
\end{definition}
When $\Gamma$ is a hyperbolic group acting on its boundary $\bgam$, every point in
$\bgam$ is a conical limit point.
The first step in our construction is to use this property to set up a
pair of good covers of $\bgam$. Then we use these covers to create a
finite-state automaton that accepts words in (a finite index subgroup
of) $\Gamma$ which code points. To produce the covers, we use the
following general lemma. Here and in what follows, metric notions like
diameter always refer to the fixed metric $\dvis$ on $\bgam$.

\begin{lemma}[Expanded neighborhoods]\label{lem:goodneighborhoods}
  For any positive $\epsilon < \frac{D}{4}$ and any $z \in \bgam$,
  there is an $\alpha_z\in \Gamma$ and a pair of open neighborhoods
  $\hat{V}_z\subset W_z$ of $z$ so that
  \begin{enumerate}
  \item\label{lemitm:smalldiam} $\diam(W_z) \le \epsilon$;
  \item\label{lemitm:bigW} $\diam(\alpha_z^{-1}W_z)>3\epsilon$; and
  \item\label{lemitm:deepnestU} $\overline{N_{1.5\epsilon}(\alpha_z^{-1}\hat{V}_z)} \subset \alpha_z^{-1}W_z$.
  \end{enumerate}
\end{lemma}
\begin{proof}
  We choose some $\epsilon<\frac{D}{4}$ where $D$ is the constant from Lemma~\ref{lem:distantpairs}.
  
  Let $z\in \bgam$.  Since $z$ is a conical limit point, we can find distinct points $a,b$ and a sequence of group elements $\{g_i\}_{i\in\bN}$ so that $g_iz\to b$ and $g_i x\to a$ uniformly away from $z$.  Up to post-composing all $g_i$ with a fixed element $g$ as in Lemma~\ref{lem:distantpairs} if necessary, we may assume $\dvis(a,b) \ge D$.  Also, since $\bgam$ is perfect, there is some point $a'\ne a$ with $\dvis(a,a') = \epsilon' <\epsilon$.  
  
  Let $W_z=B_{\epsilon/2}(z)$, so Property~\eqref{lemitm:smalldiam} is satisfied.  Let $K_z$ be the complement of $W_z$ in $\bgam$.  The set $K_z$ is compact and does not contain $z$, so for $i$ sufficiently large, we have $g_i K_z\subset B_{\epsilon'}(a)$ and $g_i z\in B_\epsilon(b)$.

  Fixing some such $i$, set $\alpha_z= g_i^{-1}$, and let
  $\hat{V}_z = \alpha_z(B_\epsilon(b))$.  Note that $B_\epsilon(a)$ contains
  $\alpha_z^{-1} K_z = \bgam - \alpha_z^{-1} W_z$.

  Since $B_\epsilon(a)$ is
  disjoint from $B_\epsilon(b)$, we have
  \[ \hat{V}_z = \alpha_z(B_\epsilon(b))\quad \subset\quad  \bgam - \alpha_z (K_z) = W_z. \]

  The set $\alpha_z^{-1}W_z = \bgam \minus \alpha_z^{-1}K_z$ contains both $b$ and $a'$, so $\diam(\alpha_z^{-1}W_z) \ge \dvis(b,a') \ge D-\epsilon > 3\epsilon$, establishing Property~\eqref{lemitm:bigW}.

 Finally, since $\dvis(b,\alpha_z^{-1}K_z)\ge D-\epsilon > 3\epsilon$, we have
  \[ \overline{N_{1.5 \epsilon}(\alpha_z^{-1}\hat{V}_z)} = \overline{B_{2.5 \epsilon}(b)} \subset (\bgam \minus \alpha_z^{-1} K_z) = \alpha_z^{-1}W_z, \]
  establishing Property~\eqref{lemitm:deepnestU}.
\end{proof}
\begin{remark}\label{rem:cell}
  The group $\Gamma$ is not required to act quasi-conformally with respect to the metric $\dgam$.  This means that in case $\bgam$ is a topological $n$--sphere, we can take $\dgam$ to be a round metric.  The sets $W_z$ can therefore be taken to be open $n$--cells whose closures are closed $n$--cells.  This will be important in the proof of Corollary~\ref{cor:cell}.
\end{remark}

\begin{definition}[Fixing $\epsilon$]\label{def:epsilon}
  Here we will fix a scale $\epsilon$ for the rest of the paper.  In order to do so we first fix, as in the statement of Theorem~\ref{thm:main}, a neighborhood $\mc{U}$ of the identity in the space of continuous self-maps of $\bgam$.  Now we fix a scale $\epsilon$ so that all of the following hold.
  \begin{enumerate}
  \item $0 < \epsilon < D/4$, where $D$ is the constant furnished by Lemma~\ref{lem:distantpairs}. 
  \item $\mathcal{U}$ contains all maps $f$ such that $d(x, f(x)) \leq \epsilon$ for all $x \in \bgam$.  
  \end{enumerate}
\end{definition}

\begin{definition}[Fixing a pair of covers]\label{def:covers}
  Since $\epsilon< D/4$, we can apply
  Lemma~\ref{lem:goodneighborhoods}.  For each $z\in \bgam$ choose a
  pair of open neighborhoods $\hat{V}_z\subset W_z$ of $z$ as in the conclusion
  of Lemma~\ref{lem:goodneighborhoods}.  Let $I \subset \bgam$ be a
  finite collection so that the sets $\{\hat{V}_z\}_{z\in I}$ cover $\bgam$.
  Parts~\eqref{lemitm:smalldiam} and~\eqref{lemitm:bigW} of Lemma~\ref{lem:goodneighborhoods} imply that the $W_z$ satisfy the following conditions.
  \begin{enumerate}[label = (C\arabic*)]
  \item\label{itm:smalldiam} $\diam(W_z) \le \epsilon$;
  \item\label{itm:bigW} $\diam(\alpha_z^{-1}W_z)>3\epsilon$ 
  \end{enumerate}
  Our sets $\hat{V}_z$ could have the inconvenient property that for some pair $y,z$, the intersection of \emph{closures} $\overline{\hat{V}}_z\cap \alpha_y^{-1}\overline{\hat{V}}_y$ is non-empty, even when the corresponding intersection $\hat{V}_z\cap \alpha_y^{-1}\hat{V}_y$  is empty, meaning that this intersection may not persist under small perturbations of $\alpha_y$.  
  However, since there are only finitely many such pairs to consider, we may modify the sets $\hat{V}_z$ slightly to change each such unstable intersection to a stable one.  Precisely, we replace each set $\hat{V}_z$ by a larger set $V_z$ so that no unstable intersections occur, and we take the enlargements small enough so that  $V_z \subset W_z$ still holds, and there are no further intersections of these sets and their images under the $\alpha_z$.  For a sufficiently small modification, condition~\eqref{lemitm:deepnestU} of Lemma \ref{lem:goodneighborhoods} will still hold if $1.5\epsilon$ is replaced by $\epsilon$.  In summary, the sets $V_z$ have the following properties:  
  \begin{enumerate}[resume,label = (C\arabic*)]
    \item\label{itm:deepnest} $\overline{N_{\epsilon}(\alpha_z^{-1}V_z)} \subset \alpha_z^{-1}W_z$.
    \item\label{itm:cleanoverlap} For any pair $y, z \in I$, the intersection
  $\overline{V}_z \cap \alpha_y^{-1}\overline{V}_y$ is either empty or
  has nonempty interior.
\end{enumerate}
\end{definition}
\putinbox{For the rest of the paper we fix the scale $\epsilon$, the indexing set $I$ and pair of covers $V_z\subset W_z$ from Definition~\ref{def:covers}.
  For each $z \in I$ we also fix the group element $\alpha_z$ from
  Lemma~\ref{lem:goodneighborhoods}.
  
}
\begin{lemma}\label{lem:axiom2}
  For each $z, y\in I$, if $\overline{V}_z$ meets  $\alpha_y^{-1}\overline{V}_y$, then $\overline{W_z}$ is contained in $\alpha_y^{-1}(W_y)$.
\end{lemma}
\begin{proof}
  Suppose $\overline{V}_z\cap \alpha_y^{-1}\overline{V}_y$ is
  non-empty.  Property~\ref{itm:deepnest} of Definition~\ref{def:covers} implies that
  $\overline{V}_z\subset W_z$, so we deduce that
  $W_z\cap \alpha_y^{-1}\overline{V}_y$ is nonempty.
  Property~\ref{itm:smalldiam} of  Definition~\ref{def:covers}
  implies $\diam(W_z) \le \epsilon$.  Thus 
  $\overline{W}_z$ is contained in $\overline{N_\epsilon(a_y^{-1}V_y)}$, which
  is contained in $\alpha_y^{-1}W_y$ by Property~\ref{itm:deepnest} of  Definition~\ref{def:covers}.  
\end{proof}

The index set $I$ can be given the structure of a directed graph as follows.
\begin{definition}[The associated automaton]
We let $\mc{G}$ be the graph with vertex set $I$, with an edge from $y$ to $z$
if and only if $\alpha_y^{-1}(\overline{V}_y)\cap \overline{V}_z \ne \emptyset$ (see Figure~\ref{fig:edgecondition}).
\begin{figure}[htbp]
  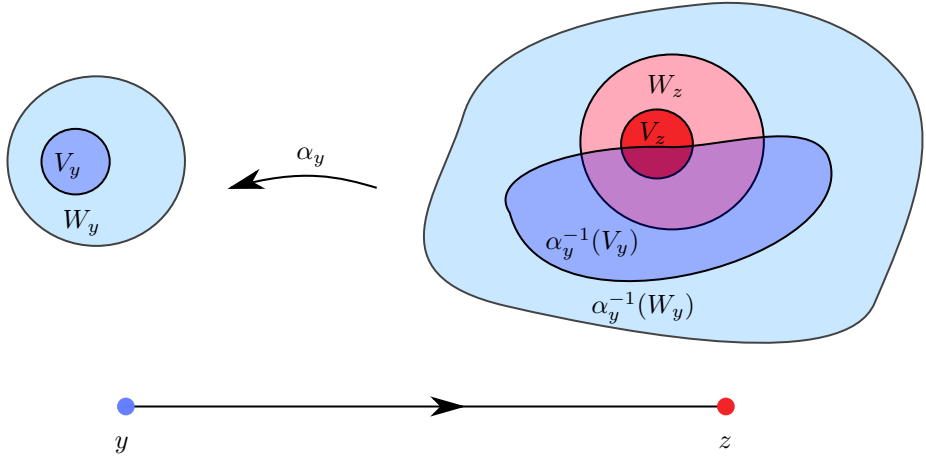
  \caption{There is an edge from $y$ to $z$ when $\overline{\alpha_y^{-1}(V_y)}$ meets $\overline{V}_z$.}
  \label{fig:edgecondition}
\end{figure}
A \emph{$\mc{G}$--coding} is an infinite sequence $\{z(k)\}_{k\in\bN}$ of points of $I$ so that there is an edge from $z(k)$ to $z(k+1)$ for each $k$.
  If 
  \begin{equation*}
    p \in \bigcap_{k = 0}^\infty\alpha_{z(1)}\cdots\alpha_{z(k)}\overline{W}_{z(k+1)}
  \end{equation*}
  we say that $\{z(k)\}_{k\in\bN}$ is a \emph{$\mc{G}$--coding of $p$}.  If $\mc{G}$ is understood, we may omit it, and speak of a \emph{coding of $p$}.
\end{definition}
\begin{remark}
Lemma~\ref{lem:axiom2} implies that for any directed path $\{z(k)\}_{k\in\bN}$, we have $\alpha_{z(k)}\overline{W}_{z(k+1)} \subset W_{z(k)}$.  Thus, the intersection 
\[ \bigcap_{k = 0}^\infty\alpha_{z(1)}\cdots\alpha_{z(k)}\overline{W}_{z(k+1)} \]
is an intersection of nested closed sets, and so by compactness of $\bgam$, it is always nonempty.  In particular every $\mc{G}$--coding is a $\mc{G}$--coding of at least one point.  We will see below that it is a $\mc{G}$--coding for \emph{only} one point.
\end{remark}
\begin{lemma}\label{lem:codingsexist}
  Every $p\in \bgam$ has a coding.
\end{lemma}
\begin{proof}
  Let $p \in \bgam$ be given.  Take $z(1)\in I$ so that $p\in \overline{V}_{z(1)}$.  Suppose that $z(1),\ldots,z(k)$ have been defined, and let $g_k = \alpha_{z(1)}\cdots \alpha_{z(k)}$.  Then there is some $z(k+1)$ so that $g_k^{-1}p\in \overline{V}_{z(k+1)}$.  Inductively we have $g_{k-1}^{-1}(p)\in \overline{V}_{z(k)}$, so $g_k^{-1}p = \alpha_{z(k)}^{-1}g_{k-1}^{-1}(p)$ is a point in the intersection of $\alpha_{z(k)}^{-1}\overline{V}_{z(k)}$ with $\overline{V}_{z(k+1)}$.  In particular there is an edge joining $z(k)$ to $z(k+1)$.  From the construction we have
  \[ p \in \bigcap_{k = 0}^\infty g_k\overline{V}_{z(k+1)} \subset \bigcap_{k = 0}^\infty g_k \overline{W}_{z(k+1)},\]
  so the sequence $\{z(k)\}_{k\in\bN}$ is a $\mc{G}$--coding of
  $p$.
\end{proof}
\begin{lemma}[Bounded backtracking property] \label{lem:boundedbacktrack}
  Let $\{z(k)\}_{k\in\bN}$ be a $\mc{G}$--coding.  For any $k$ define
  \[ U_{k} = \alpha_{z(1)}\cdots\alpha_{z(k-1)}W_{z(k)}. \] Then
  $U_{k+1}$ is a proper subset of $U_k$ for any $k\ge 1$.  Moreover, in the sequence
  $\{g_k = \alpha_{z(1)}\cdots\alpha_{z(k)}\}_{k\in\bN}$, no element $g_k$
  is repeated more than $\#I$ times. 
\end{lemma}

\begin{proof}
  Fix $k\ge 1$.  By the definition of the graph $\mc{G}$, we
  have
  \[\alpha_{z(k)}^{-1}(\overline{V}_{z(k)})\cap \overline{V}_{z(k+1)}
    \ne \emptyset.\] By Lemma~\ref{lem:axiom2}, we have
  $W_{z(k+1)}\subset \alpha_{z(k)}^{-1}(W_{z(k)})$.  By
  Properties~\ref{itm:smalldiam} and~\ref{itm:bigW} of
  Definition~\ref{def:covers}, the set $W_{z(k+1)}$ has diameter at
  most $\epsilon$, whereas the set $\alpha_{z(k)}^{-1} W_{z(k)}$ has
  diameter at least $3\epsilon$.  In particular the inclusion
  $W_{z(k+1)}\subset \alpha_{z(k)}^{-1}(W_{z(k)})$ is
  proper. Multiplying on the left by
  $\alpha_{z(1)}\cdots \alpha_{z(k)}$ then gives a proper inclusion
  $U_{k+1}\subset U_k$.

  To see the last assertion, suppose $\#\{k \mid g_k = g\}>\#I$ for some $g$.  Then there must be distinct $k, k'$ so that $g_k = g_{k'} = g$ and  $W_{z(k)} = W_{z(k')}$.  For these indices, $U_k  = U_{k'} = g W_{z(k)}$, a contradiction to proper nesting.
\end{proof}

\section{Properties of $\mc{G}$--codings} \label{sec:codingproperties}
In this section we establish key properties of $\mc{G}$--codings to be
used in the proof of the main theorem.  Our first goal is to show that a sequence in $\Gamma$ defined by a coding lies a uniformly bounded Hausdorff distance from a geodesic ray based at identity.  

\subsection{Codings and quasi-geodesics} 
We begin with two general lemmas on hyperbolic spaces.  
\begin{lemma}\label{lem:existsc}
  Let $\epsilon_0 < D/4$.  There is a $c_0>0$ so that for
  any $p\in \bgam$, there are points $q,q'$ so that
  $\dvis(p,q) > \epsilon_0$ and $\dvis(p,q') > \epsilon_0$, and
  $\dvis(q,q')\ge c_0$.
\end{lemma}
\begin{proof}
  We recall that $D\le \diam(\bgam)$ is the bound from Lemma~\ref{lem:distantpairs} so that any pair of points in $\bgam$ can be translated to a pair whose distance is at least $D$.  In particular, there are points $N,S\in \bgam$ with $\dvis(N,S) \ge D$.  Since $\bgam$ is perfect, there are points $N'\in B_{\epsilon_0}(N)\minus \{N\}$ and $S'\in B_{\epsilon_0}(S)\minus\{S\}$.  Let $c_0 = \min\{\dvis(N,N'),\dvis(S,S')\}$.

  If $p$ is further than $\epsilon_0$ from both $N$ and $S$, we may take $q = N$, $q' = S$.  Otherwise we may suppose after relabeling that $\dvis(p,N) \le \epsilon_0$.  But then
  \[\dvis(p,S) \ge D-\epsilon_0 \ge 3\epsilon_0\] and \[\dvis(p,S') \ge D - 2\epsilon_0 \ge  2\epsilon_0,\] so we may take $q = S$ and $q' = S'$.
\end{proof}

\begin{lemma}\label{lem:tripodnature}
  Given $\delta \ge 0$, there is a constant $T$ depending only on
  $\delta$ so that the following holds.  Let $(X, d_X)$ be a
  $\delta$--hyperbolic metric space.  Let $p,q,r\in \partial X$, and for each
  pair of distinct points $x,y \in \{p,q,r\}$, let $(x, y)$ be a
  bi-infinite geodesic joining $x$ to $y$.

  Then for each $x, y \in \{p, q, r\}$ distinct, there are points
  $c_{xy} \in (x, y)$ such that
  \begin{enumerate}
  \item\label{itm:close} $\diam(\{c_{pq}, c_{qr}, c_{pr}\})\le T$;
  \item\label{itm:far}  If $w$ lies in the sub-ray $[c_{pq},p)\subset (p,q)$, then
    \[ d(w, (q,r) ) \ge d(w, c_{pq}) - T\]
    and similar statements hold with $p,q,r$ permuted.
  \end{enumerate}
\end{lemma}

\begin{figure}[htbp]
  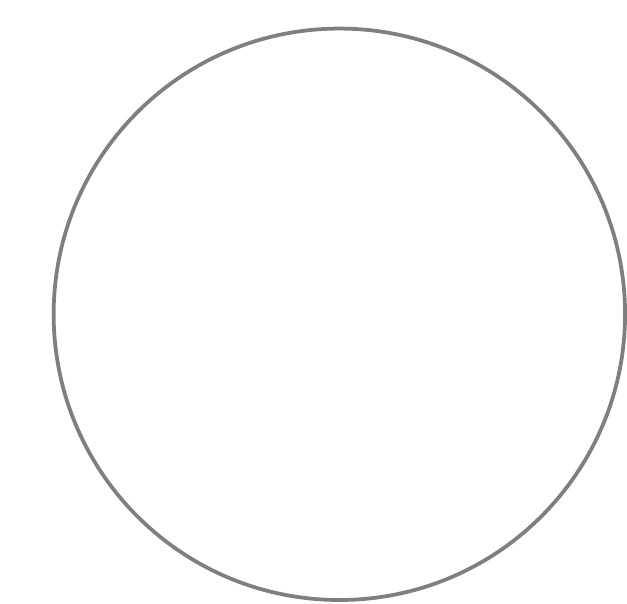
  \caption{Lemma~\ref{lem:tripodnature}.}
  \label{fig:tripodnature}
\end{figure}
\begin{proof}
  We may assume that for any $x, y$ distinct in $\{p, q, r\}$, we have
  $(x, y) = (y, x)$.  Let $c$ be any point which lies within $2\delta$
  of all three geodesics, and for each $\{x, y\} \subset \{p, q, r\}$,
  choose a closest point $c_{yx} = c_{xy} \in (x,y)$ to $c$.  These
  points clearly satisfy \eqref{itm:close}.

  The point
  $c_{xy}$ cuts $(x, y)$ into rays $[c_{xy}, x)$ and $[c_{xy},
  y)$. For each $z\in \{p,q,r\}$, choose a geodesic ray $[c, z)$ from
  $c$ to $z$.  Let $Y$ be the tripod which is the union of these three rays, and let $\Delta$ be the union of the geodesics $(p,q),(q,r)$ and $(p,r)$.   We define
  a map $\pi\from \Delta \to Y$ so that
    $\pi(\{c_{pq}, c_{qr}, c_{rp}\}) = z$, and so that $\pi$ sends the
  ray $[c_{xy}, y)$ isometrically to the ray $[c, y)$. It is
  straightforward to see that $d_X(\pi(w), w)\le 4\delta$ for any
  $w \in \Delta$, so $\pi$ is a $(1,8\delta)$--quasi-isometry.
  
  Let $\tilde Y$ be an abstract infinite tripod, the wedge of three rays.  We
  claim the obvious map from $\tilde Y$ to $Y$, isometric on each ray, is a
  $(1,8\delta)$--quasi-isometric embedding.  Indeed, the embedding is
  clearly distance non-increasing, and it is enough to consider points
  on distinct rays.  Up to relabeling, we may assume that these two
  rays are $[c, p)$ and $[c, q)$. We consider unit-speed
  parameterizations $\tau_p$ and $\tau_q$ of these rays, and
  unit-speed parameterizations $\alpha_p$, $\alpha_q$ of the rays
  $[c_{pq}, p)$ and $[c_{pq}, q)$. Then, for any
  $s, t \ge 0$, we have
  $d_X(\alpha_p(s), \tau_p(s)) \le 4\delta$ and
  $d_X(\alpha_q(t), \tau_q(t)) \le 4\delta$. Since $(p, q)$ is
  geodesic, $d_X(\tau_p(s),\tau_q(t))\ge s + t - 8\delta$.

  Combining the maps from the last two paragraphs, we see that the triangle $\Delta$ is $(1,16\delta)$--quasi-isometric
  to the abstract tripod $\tilde Y$.  Since \eqref{itm:far} holds in the tripod with $T=0$ (taking all the $c_{xy}$ to be the central point of the tripod), it holds in $\Delta$ with $T = 32\delta$.
\end{proof}

\begin{definition}
  The points $c_{pq}, c_{qr}, c_{pr}$ from
  Lemma~\ref{lem:tripodnature} will be referred to as \emph{central
    points} of the ideal triangle with vertices $p$, $q$, $r$.
\end{definition}

The next lemma 
is inspired by Proposition 5.11 from \cite{Weisman}.
\begin{lemma}\label{lem:quasigeodesic}
Let $X$ denote the Cayley graph of $G$ with generating set $\mathcal{S}$.  Fix a $\mc{G}$--coding $\{z(k)\}_{k\in\bN}$.  For ease of notation, let $\alpha_k$ denote $\alpha_{z(k)}$, and let $W_k$ denote $W_{z(k)}$.   
 Then the set $\{g_k = \alpha_{1}\cdots\alpha_{k}\}_{k\in\bN}$ is (uniformly) finite Hausdorff distance from a geodesic ray based at the identity in $X$.  If $\{z(k)\}_{k\in\bN}$ is a $\mc{G}$--coding for $p$, then this geodesic ray tends to $p$. 
\end{lemma}

\begin{proof}
  Let $\{z(k)\}_{k\in\bN}$ be a $\mc{G}$--coding for $p$ and let
  $g_k = \alpha_{1}\cdots\alpha_{k}$. The distances
  $d_X(g_k, g_{k+1})$ are uniformly bounded, so $\{g_k\}_{k\in\bN}$ is
  uniformly finite Hausdorff distance from the image of some path in
  $X$.  Lemma~\ref{lem:boundedbacktrack} implies this path must be a
  proper ray.  So it suffices to show that the sequence
  $\{g_k\}_{k \in \bN}$ lies within a uniformly bounded neighborhood
  of a geodesic ray based at the identity in $X$, tending towards $p$.

The general strategy of proof is as follows: we first show that
$\{g_k\}$ is contained in a uniformly bounded neighborhood of a sub-ray
of a geodesic between $p$ and some point $q \in \bgam$, and that this
sub-ray also passes close to the identity vertex in $X$.  By varying
the point $q$ and running the same argument, we show that this sub-ray
does not extend too far in the direction of $q$ and conclude that it
is close to a ray based at the identity vertex (see Figure
\ref{fig:close_ray} for a schematic).

Let $c_1$ be a positive lower bound for the distances
  $\dvis(\bgam \minus \alpha_x^{-1}W_x,\overline{W}_y)$ as $(x,y)$
  ranges over the directed edges of $\mc{G}$.  Let $c_0$ be the constant from Lemma~\ref{lem:existsc} (applied with $\epsilon_0=\epsilon$) and let $c = \min\{c_0,c_1,\epsilon\}$.  
Take $C > 0$ large enough 
 (depending only on $c$) so that for any pair of points $x,y$ with
  $\dvis(x,y)\ge c$, any geodesic joining $x$ to $y$ in $X$ comes
  within $C$ of the identity vertex $\idgamma$.  The existence of such a $C$ is immediate if $\dvis$ is a visual metric based at $\idgamma$; since the boundary is compact, such a $C$ will exist for any metric.

  By Lemma~\ref{lem:existsc}, there are points $q,q'$ so that
  $\dvis(p,q)>\epsilon$, $\dvis(p,q')>\epsilon$, and
  $\dvis(q,q')>c_0$.  In particular, any geodesic joining a pair of
  distinct points in $\{p,q,q'\}$ meets the ball of radius $C$ around
  the identity.
  
 For each $k \in \bN$, let $\overline{U}_k = g_{k-1}\overline{W}_{k}$.  This
 gives a family of nested closed neighborhoods of $p$, each of which has diameter $\le \epsilon$.
  Since $\dvis(p,q)$ and $\dvis(p,q')$ are strictly larger than $\epsilon$, neither $q$ nor $q'$ is contained in any of the  sets $\overline{U}_k$.  In particular,
  \begin{equation}
    \label{eq:notinW1}
    q,q' \notin \overline U_k \subset \overline W_1.
  \end{equation}

  Now let $z$ be either $q$ or $q'$.  For any $k$ we have
  $g_k^{-1}(p) \in g_k^{-1}\overline{U}_{k+1} = \overline{W}_{k+1}$
  and, by~\eqref{eq:notinW1},
  $g_k^{-1}(z) \notin g_k^{-1}\overline{U}_{k} =
  \alpha_{k}^{-1}(\overline{W}_{k})$.  Thus
  $\dvis(g_k^{-1}(p),g_k^{-1}(z)) \ge c_1 \ge c$, which by our choice of $C$ means that
  any geodesic from $g_k^{-1}(p)$ to $g_k^{-1}(z)$ passes within distance $C$ of $e$; in other words,
if we fix a geodesic
  $\gamma$ from $z$ to $p$, we have
  $\dgam(\idgamma, g_k^{-1}\gamma)\le C$ for each $k$. Multiplying by
  $g_k$ we see that $\dgam(g_k, \gamma) \le C$ for each $k$.  In
  particular the sequence $\{g_k\}_{k \in \bN}$ is contained in the closed
  $C$--neighborhood of $\gamma$.

  Now fix bi-infinite geodesics $(q, q')$, $(p, q)$ and $(p, q')$, and
  consider central points $c_{pq}, c_{pq'}$ and $c_{qq'}$ for the
  ideal triangle with vertices $p, q, q'$, given by
  Lemma~\ref{lem:tripodnature}. We claim that these central points lie
  uniformly close to $\idgamma$. To see this, for each
  $\{x, y\} \subset \{p, q, q'\}$, we let $w_{xy}$ be a point lying on
  the geodesic from $x$ to $y$ such that $w_{xy}$ is within distance
  $C$ of the identity. If $\{x, y, z\} = \{p, q, r\}$, the geodesic
  segment $[c_{xy}, w_{xy}]$ lies in either the sub-ray $[c_{xy}, x)$
  or $[c_{xy}, y)$; up to relabeling, we assume it lies in
  $[c_{xy}, x)$. Then, by Lemma~\ref{lem:tripodnature}, we have
  \[
    2C \ge \dgam(w_{xy}, w_{yz}) \geq \dgam(w_{xy}, (y, z)) \geq
    \dgam(w_{xy}, c_{xy}) - T.
  \]
  This gives the estimate $\dgam(c_{xy}, \idgamma) < 3C + T$, as desired.  

  Using this bound, we conclude that the rays
  $[c_{pq}, p) \subset (q, p)$ and $[c_{pq'}, p) \subset (q', p)$ are
  both Hausdorff distance at most $3C + T + 2\delta$ from any geodesic
  ray $[\idgamma, p)$ from $\idgamma$ to $p$.

  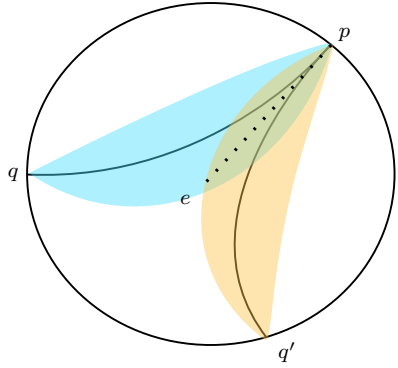
\begin{figure*}
   \labellist 
  \footnotesize \hair 2pt
     \pinlabel $p$ at 255 250
\pinlabel $q$ at -5 140
    \pinlabel $q'$ at 210 0
     \pinlabel $\idgamma$ at 130 120
   \endlabellist
     \centerline{ \mbox{
 \includegraphics[width = 2in]{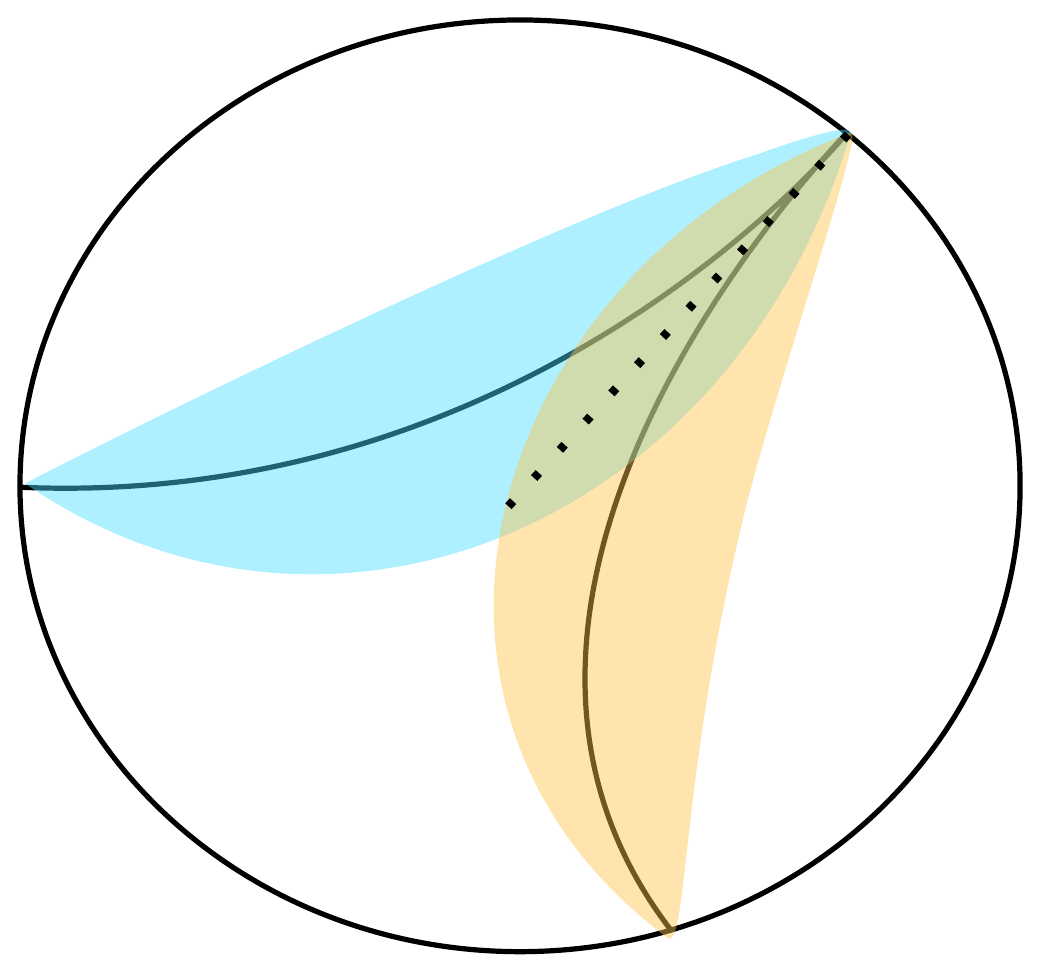}}}
 \caption{If $\{g_k\}$ lies in a bounded neighborhood of both $pq$ and $pq'$, then it must be near a ray to $p$ based at $\idgamma$}
  \label{fig:close_ray}
  \end{figure*}

  Let $z\in \{g_k\}_{k\in\bN}$.  Then $z$ lies within $C$ of both $(q, p)$ and
  $(q', p)$. If $z$ lies within $C$ of $[c_{pq}, p) \subset (q, p)$,
  then it lies in the $4C + T + 2\delta$--neighborhood of
  $[\idgamma, p)$.  Otherwise, $z$ lies within distance $C$ of a point
  $w$ on $[c_{pq}, q)$. Then by Lemma~\ref{lem:tripodnature} again, we
  have
  \[
    \dgam(z, c_{pq}) \le \dgam(w, c_{pq}) + C \le \dgam(w, (p, q')) +
    C + T \le \dgam(z, (p, q')) + 2C + T.
  \]
  Since $\dgam(c_{pq}, \idgamma) \le 3C + T$ and
  $\dgam(z, (p, q')) \le C$, we have $\dgam(z, \idgamma) \le 6C +
  2T$. In either case $z$ lies in the
  $6C + 2T + 2\delta$--neighborhood of the ray $[\idgamma, p)$.

\end{proof}
\begin{corollary} \label{cor:unique} 
  For any $\mc{G}$--coding $\{z(k)\}_{k\in\bN}$, the intersection 
 \[ \bigcap_{k = 0}^\infty \alpha_{z(1)}\cdots \alpha_{z(k)}\overline{W}_{z(k+1)} \]
 is a singleton.     
\end{corollary}
\begin{proof}
  By Lemma~\ref{lem:quasigeodesic}, $\{z(k)\}_{k\in\bN}$ determines a set $\{g_k\}_{k\in\bN}$ which is finite Hausdorff distance from a geodesic ray with a unique endpoint $p$.   If $q$ is a point in the intersection, then $\{z(k)\}_{k\in\bN}$ codes $q$, and therefore by Lemma~\ref{lem:quasigeodesic}, $q=p$.
\end{proof}

Combining Lemmas~\ref{lem:quasigeodesic} and~\ref{lem:codingsexist} we have the following.
\begin{corollary}
  The set $\{\alpha_z\}_{z\in I}$ generates a finite index subgroup of
  $\Gamma$.
\end{corollary}

\subsection{Uniform contraction} 
This section contains the primary technical application of our work above.  
The following is an immediate consequence of Lemma \ref{lem:quasigeodesic}.
\begin{lemma}\label{lem:F}
  There is a finite set $F\subset G$ so that the following holds.  Let $p\in \bgam$, and let $\{z(k)\}_{k\in\bN}, \{y(k)\}_{k\in\bN}$ be two $\mc{G}$--codings of $p$.  For $k>0$ let $g_k = \alpha_{z(1)}\cdots\alpha_{z(k)}$ and $h_k = \alpha_{y(1)}\cdots\alpha_{y(k)}$.  For any $k>0$ there is some $n(k)$ so that
  \begin{equation}  \label{eq:differbyF}
  h_k =  g_{n(k)}f_k \text{ for some } f_k \in F. 
  \end{equation}
\end{lemma}
Symmetrically, there exists $m(k) \in \mathbb{N}$ such that 
$g_k = h_{m(k)}f'_k \text{ for some } f'_k \in F$, but it is the equality~\eqref{eq:differbyF} that we will make use of in the next step of the proof.  

The following lemma is the technical heart of our main theorem.  

\begin{lemma}[Uniform Contraction] \label{lem:uniform} Using the
  notation from Lemma~\ref{lem:F}, there exists a uniform $N>0$ so
  that for any such pair of such sequences $\{g_k\}_{k \in \bN}$ and
  $\{h_k\}_{k \in \bN}$ and any $k>0$, we have
  \begin{equation}
    \label{eq:nestingcond}
   g_{n(k)+N} \overline{W}_{z(n(k)+N+1)} \subset h_{k} W_{y(k+1)} .
 \end{equation}
\end{lemma}

We emphasize that by ``uniform'' we mean that the constant $N$ is independent of the point being coded and of the two chosen codings of the point.  

\begin{proof}
  The proof is by contradiction.  We therefore assume we have a sequence of natural numbers $N$
  tending to $\infty$ so that for each $N$ in the sequence there is a point $p^{N}\in\bgam$ with two distinct codings $\{z^{N}(k)\}_{k\in\bN}$ and $\{y^{N}(k)\}_{k\in\bN}$, so that the
  statement~\eqref{eq:nestingcond} fails for this pair of codings and some $k=k_N \in \bN$.

Before continuing the proof, we clarify notation: 
Our sequences of codings $\{z^{N}(k)\}_{k\in\bN}$ and
$\{y^{N}(k)\}_{k\in\bN}$ are fixed and indexed by $N$.  As in Lemma~\ref{lem:F}, for 
$g^{(N)}_k = \alpha_{z^{N}(1)}\cdots\alpha_{z^{N}(k)}$ and $h^{(N)}_k = \alpha_{y^{N}(1)}\cdots\alpha_{y^{N}(k)}$
we have indices $n^{N}(k) \in \bN$ and elements $f^{(N)}_k$ of the finite set $F$ such that for each $N,k$ there is an equality
  \begin{equation}
    \label{eq:finiteset}
 h^{(N)}_{k} = g^{(N)}_{n^{N}(k)}f^{(N)}_k.
\end{equation}
\newcommand{\nN}{n_N}
To decrease notation slightly, we shorten $n^{N}(k_N)$ to $\nN$.
The failure of the nesting condition~\eqref{eq:nestingcond} can then be expressed as
  \begin{equation}
    \label{eq:notsubset}
    g^{(N)}_{\nN+N} \overline{W}_{z^{N}(\nN+N+1)}
    \not\subset h^{(N)}_{k_N} W_{y^{N}(k_N + 1)}.
  \end{equation}
  We immediately pass to a subsequence so that
  \begin{equation}\tag{f}
    \label{cond:f} \mbox{ $f^{(N)}_{k_N}$ is constant, equal to $f \in F$.  }
  \end{equation}
  We further refine our subsequence so the following three
  conditions are satisfied.
  \begin{align}
  \tag{z} & \mbox{ The sets $W_{z^{N}(\nN+N+1)}$ are constant, equal to some $W_z$.}\\
  \tag{y} & \mbox{ The sets $W_{y^{N}(k_N+1)}$ are constant, equal to some $W_y$. }\\
  \tag{$\ast$}\label{cond:W} & \mbox{ The sets $\alpha_{y^N(k_N+1)}W_{y^N(k_N+2)}$ are constant, equal to some $W$.}
  \end{align}
  This is possible because there are only finitely many possibilities
  for $W_{z^{N}(\nN+N+1)}$ and $W_{y^{N}(k_N+1)}$ (being elements
  of the cover) and also for $\alpha_{y^N(k_N+1)}W_{y^N(k_N+2)}$,
  being a translate of an element of the cover by one of finitely many
  elements.  A key property we will use at the end of the proof is that
  \begin{equation}\label{eq:WinWy}
    \overline{W}\subset W_y.
  \end{equation}
  
  For these values of $N$, we can multiply each side of~\eqref{eq:notsubset} on the left by $(g^{(N)}_{\nN})^{-1}$ and use \eqref{eq:finiteset} together with condition \eqref{cond:f} to obtain
  \begin{equation}
    \label{eq:simple}
 \alpha_{z^{N}(\nN+1)}\cdots\alpha_{z^{N}(\nN + N)} \overline{W}_z \not\subset f W_y.
  \end{equation}
  
  For each $N$, consider the sub-coding from $\{z^{N}(k)\}_{k\in\bN}$
  formed by terms $\nN + 1$ through infinity.  We denote this coding
  by $\{\gamma^N(k)\}_{k\in\bN}$.  In other words, we define
  $\gamma^N(k) = z^{N}(\nN + k)$.  After passing to a subsequence
  $\{N(j)\}_{j\in\bN}$ one final time, we may obtain a sequence of
  codings $\{\gamma^{N(j)}(k)\}_{k\in\bN}$ so that for all $l \geq j$,
  the initial segment $\{\gamma^{N(l)}(1),\ldots,\gamma^{N(l)}(j)\}$
  is independent of $l$, so equal to the first $j$ terms of
  $\{\gamma^{N(j)}(k)\}_{k\in\bN}$.  That is, the codings
  $\{\gamma^{N(j)}(k)\}_{k\in\bN}$ converge to a coding
  $\{\gamma^\infty(k)\}_{k\in\bN}$, which codes a unique point
  $p^{\infty}\in \bgam$.
  
For our subsequence $N(j)$, the non-containment in \eqref{eq:simple} takes the slightly simpler form
\begin{equation*}
 \alpha_{\gamma^{N(j)}(1)}\cdots\alpha_{\gamma^{N(j)}(N(j))} \overline{W}_z \not\subset f W_y,
\end{equation*}
which we can rewrite as
\begin{equation*}
\left( \alpha_{\gamma^{\infty}(1)}\cdots\alpha_{\gamma^\infty(j)}\right)\left(\alpha_{\gamma^{N(j)}(j+1)}\cdots \alpha_{\gamma^{N(j)}(N(j))}\right) \overline{W}_z \not\subset f W_y.
\end{equation*}
Since $\{\gamma^{N(j)}(k)\}_{k\in\bN}$ is a coding, 
\begin{equation*}
  \left(\alpha_{\gamma^{\infty}(1)}\cdots\alpha_{\gamma^\infty(j)}\right)\left(\alpha_{\gamma^{N(j)}(j+1)}\cdots \alpha_{\gamma^{N(j)}(N(j))} \right)\overline{W}_z\subset  \alpha_{\gamma^{\infty}(1)}\cdots\alpha_{\gamma^\infty(j-1)}W_{\gamma^\infty(j)}, 
\end{equation*}
so we must therefore have 
\begin{equation}\label{eq:pinftynbhds}
\alpha_{\gamma^{\infty}(1)}\cdots\alpha_{\gamma^\infty(j-1)}\overline{W}_{\gamma^\infty(j)}\not\subset f W_y.
\end{equation}
Since $\{\gamma^\infty(k)\}_{k\in\bN}$ is a coding for $p^{\infty}$, the sets on the left hand side of~\eqref{eq:pinftynbhds} give a nested basis of closed neighborhoods of $p^{\infty}$ and we must have
\begin{equation}
  \label{eq:notinrhs}
  p^{\infty}\not\in f W_y.
\end{equation}

On the other hand, since $\{z^{N}(k)\}_{k\in\bN}$ is a $\mc{G}$--coding of
$p^{N}$, for every $j$ we have
\begin{equation*}
  p^{N(j)} \in \alpha_{z^{N(j)}(1)} \cdots \alpha_{z^{N(j)}(n_{N(j)}
    + j - 1)} W_{z^{N(j)}(n_{N(j)} + j)},
\end{equation*}
or equivalently (multiplying both sides on the left by $(g^{(N(j))}_{n_{N(j)}})^{-1}$):
\begin{align*}
  (g^{(N(j))}_{n_{N(j)}})^{-1} p^{N(j)}
  &\in \alpha_{z^{N(j)}(n_{N(j)} + 1)} \cdots \alpha_{z^{N(j)}(n_{N(j)} +
    j - 1)} W_{z^{N(j)}(n_{N(j)} + j)}\\
  &= \alpha_{\gamma^\infty(1)} \cdots \alpha_{\gamma^\infty(j - 1)} W_{\gamma^\infty(j)}.
\end{align*}
Again, since $\{\gamma^\infty(k)\}_{k\in\bN}$ codes $p^\infty$, this last
sequence of sets gives a nested neighborhood basis for $p^\infty$ and
we must have
\begin{equation*}
  \lim_{j\to\infty} (g^{(N(j))}_{n_{N(j)}})^{-1} p^{N(j)} = p^{\infty}.
\end{equation*}

We also know that
\begin{equation}\label{eq:gnpn}
  (g^{(N)}_{\nN})^{-1} p^{N}\in (g^{(N)}_{\nN})^{-1} h^{(N)}_{k_N+1}W_{y^N(k_N + 2)}
\end{equation}
for any $N$, since $ h^{(N)}_{k_N+1}W_{y^N(k_N + 2)}$ is a
 neighborhood of $p^{N}$.  By our assumptions
\eqref{cond:f} and \eqref{cond:W} on our chosen subsequence, the
right-hand side of~\eqref{eq:gnpn} is always equal to a constant
$fW$. But we have just seen that a subsequence of the left-hand side
converges to $p_\infty$, so we must have $p^{\infty}\in f\overline{W}$.  Because of~\eqref{eq:WinWy} this implies
\begin{equation*}
  p^{\infty}\in fW_y,
\end{equation*}
contradicting~\eqref{eq:notinrhs}.

\end{proof}

In our application of the uniform contraction principle above, we will use the following key observation: 

\begin{remark}[A finite list of nesting conditions suffices]
  \label{rem:finitenestingcond}
  The conditions appearing in~\eqref{eq:nestingcond} (as we vary over all possible points and codings) appear to be
  infinite in number, but multiplying both sides of~\eqref{eq:nestingcond} on the left by
  $g_{n(k)}^{-1}$ gives an equivalent condition of the form
  \begin{equation}
    \label{eq:nestingcond2}
    \alpha_{z(n(k)+1)}\cdots\alpha_{z(n(k)+N)} \overline{W}_{z(n(k)+N+1)} \subset f_kW_{y(k+1)}.
  \end{equation}
  Since $N$ is bounded, there are only finitely many possible conditions of this form.  
\end{remark}

\subsection{Conjugating by generators} 
For each generator $s \in \mathcal{S}$, consider the action $\rho_s$
of $\Gamma$ on $\bgam$ given by conjugating the standard action by
$s$.  The sets $sW_z$ and elements $s \alpha_z s^{-1}$ (for $z \in I$)
give a graph $\mc{G}^s$ that is naturally isomorphic to $\mc{G}$ and
codes points for the conjugated action.  Under the natural
isomorphism, a coding of $sp$ in $\mc{G}^s$ corresponds to a coding of
$p$ in $\mc{G}$.  Furthermore, if $\{ z(k) \}_{k\in\bN}$ is a $\mc{G}^s$--coding of $sp$,
the corresponding path
$\{sg_ks^{-1} = s\alpha_{z(1)}s^{-1} \cdots s\alpha_{z(k)}s^{-1}\}_{k\in\bN}$ in
$\Gamma$ is a (uniformly, depending only on $s$) bounded distance from
a geodesic ray in $\Gamma$ tending to $sp$, since it is the image of the path
$\{g_k = \alpha_{z(1)} \cdots \alpha_{z(k)}\}_{k\in\bN}$ under conjugacy by $s$.

From the above uniformity we obtain the following analogue of Lemma~\ref{lem:F}.
\begin{lemma}\label{lem:Fs}
  For each $s\in S$ there is a finite set $F_s\subset G$ so that the following holds, for any $p\in \bgam$.  Let $\{z(k)\}_{k\in\bN}$ be a $\mc{G}$--coding of $sp\in \bgam$, and let $\{y(k)\}_{k\in\bN}$ be a $\mc{G}$--coding of $p$ (equivalently $\{y(k)\}_{k\in\bN}$ is a $\mc{G}^s$--coding of $sp$).
  For $k>0$ let
  \begin{equation*} 
g_k = \alpha_{z(1)}\cdots \alpha_{z(k)} \text{ and } h_k = s\alpha_{y(1)}s^{-1} \cdots s\alpha_{y(k)}s^{-1}.
\end{equation*}
For any $k>0$ there is some $n(k)$ so that
  \begin{equation*}  
  h_k =  g_{n(k)}f_k \text{ for some } f_k \in F_s. 
  \end{equation*}
\end{lemma}

Analogous to Lemma~\ref{lem:uniform} we have the following.
\begin{corollary} \label{cor:uniform_for_conj} Using the notation from Lemma~\ref{lem:Fs}, there exists a uniform $N_s>0$ so that for any pair of such sequences $\{g_k\}_{k\in\bN}$ and $\{h_k\}_{k\in\bN}$ and any $k>0$, we have
  \begin{equation} 
  \label{eq:nestingcond_s}
   g_{n(k)+N_s} \overline{W}_{z(n(k)+N_s+1)} \subset h_{k} sW_{y(k+1)} .
 \end{equation}
\end{corollary} 
  The proof follows that of Lemma~\ref{lem:uniform} almost verbatim.  The main change is that the constant sets $W_y$ and $W$ from that proof are defined slightly differently as
  $W_y = sW_{y^N(k_N+1)}$ and $W = s\alpha_{y^N(k_N + 1)}W_{y^N(k_N+2)}$.

\begin{remark} \label{rem:nesting_s}
As in Remark \ref{rem:finitenestingcond}, the equations of \eqref{eq:nestingcond_s} reduce to only finitely many conditions; the fact that  
$h_k =  g_{n(k)}f_k \text{ for some } f_k \in F_s$ means that each condition is equivalent to one of the form 
\begin{equation} \label{eq:finite_s}
 \alpha_{z(n(k)+1)} \ldots  \alpha_{z(n(k)+N_s)}\overline{W}_{z}  \subset fs W_{y} 
 \end{equation} 
for some $f \in F_s$, and $W_z, W_y$ in our finite set, and $z(n(k)+1), \ldots z(n(k)+N_s)$ a path of length $N_s$ in $\mc{G}$.  
\end{remark} 

\subsection{Relationship with annulus systems} 
We conclude this section by indicating how Bowditch's framework of
\emph{annulus systems} on spaces with a convergence group action can give
an alternative strategy towards the proof of the uniform contraction lemma
(Lemma~\ref{lem:uniform}). In fact, one could use Bowditch's framework
to prove a stronger version of Lemma \ref{lem:quasigeodesic}, showing
that for any $\mc{G}$--coding $\{z(k)\}_{k\in\bN}$, the map $\bN \to \Gamma$
given by $k \mapsto \alpha_{z(1)} \cdots \alpha_{z(k)}$ is actually a
(uniform) quasi-geodesic embedding.  However, since we do not need this
stronger statement anywhere in the paper, and the setup
is somewhat involved, we only provide a sketch of the argument for the purpose
of describing the relationship.

The general setting is as follows.  
In \cite{Bowditch98}, Bowditch showed that whenever $\Gamma$ acts on a
perfect compact metrizable space $Z$ as uniform convergence group,
then it is possible to completely recover a word-hyperbolic metric on
$\Gamma$ from the data of the convergence group action. To do so,
Bowditch defines a notion of a \emph{system of annuli} on $Z$, and
relates the topological behavior of such a system to a
$\Gamma$--invariant \emph{cross-ratio} (which in turn defines a
hyperbolic metric on the space of triples in $Z$).

\begin{definition}
  An \emph{annulus} in $\bgam$ is an ordered pair $A = (A^-, A^+)$ of
  disjoint closed subsets of $\bgam$, such that
  $A_- \cup A_+ \ne \bgam$. A \emph{symmetric system of annuli} is a
  collection $\mc{A}$ of annuli such that if $A = (A^-, A^+)$ is in
  $\mc{A}$, then $-A = (A^+, A^-)$ is also in $\mc{A}$.

  A sequence of annuli $A_1, \ldots, A_n$ is said to be \emph{nested}
  if $A_i^-$ contains $\bgam - A_{i+1}^+$ for all $1 \le i < n$.
\end{definition}

The terminology is inspired by the example $\bgam \simeq S^2$, where
an ``annulus'' consisting of a pair of disjoint closed disks
determines an annulus (in the usual sense) in $S^2$. Any system of
annuli $\mc{A}$ defines a four-point cross-ratio on $\bgam$, denoted
$(\cdot, \cdot ; \cdot, \cdot)_{\mc{A}}$: for $x, y, z, w \in \bgam$,
the cross-ratio $(x, y ; z, w)_{\mc{A}}$ is the maximum length of a
nested sequence of annuli $A_1, \ldots, A_n$ in $\mc{A}$ such that
$A_1^+$ contains $\{x,y\}$ and $A_n^-$ contains $\{z, w\}$.

Work of Bowditch (see \cite{Bowditch98}, Proposition 4.8 and Sections
6 and 7) shows that whenever $\mc{A}$ is a symmetric annulus system
satisfying certain hypotheses, then the cross-ratio
$(\cdot, \cdot ; \cdot, \cdot)_{\mc{A}}$ is within bounded additive
and multiplicative error of a \emph{standard} cross-ratio
$(\cdot, \cdot ; \cdot, \cdot)$ on $\bgam$. This cross-ratio is
defined by realizing $\bgam$ as the Gromov boundary of a hyperbolic
metric space $X'$ (which is $\Gamma$--equivariantly quasi-isometric to
the Cayley graph $X$ of $\Gamma$) and defining $(a, b; c, d)$ to be
the minimum distance between a geodesic in $X'$ joining $a$ to $b$ and
a geodesic joining $c$ to $d$.

We can use our sets $W_y$ for $y \in I$ to define a suitable annulus
system: for each $y \in I$, take $A_y = (A_y^+, A_y^-)$, with
\[
  A_y^+ = \bgam - W_y, \qquad A_y^- = \bigcup_{y \to z} \alpha_y
  \overline{W_z}.
\]
Then take
$\mc{A} = \{\gamma A_y \mid \gamma \in \Gamma, y \in I\} \cup \{- (\gamma
A_y) \mid \gamma \in \Gamma, y \in I\}$. One can then check that this system of
annuli meets all of Bowditch's conditions, and therefore the induced
cross-ratio $(\cdot, \cdot ; \cdot, \cdot)_{\mc{A}}$ approximates
$(\cdot, \cdot ; \cdot, \cdot)$.

Now let $\{z(k)\}_{k \in \bN}$ be a $\mc{G}$--coding, and for each
$n \in \bN$, let $q_n, q_n'$ be a pair of distinct points in
$U_{z(n + 1)}$, joined by a geodesic passing within some uniformly
bounded distance of the identity. Let
$g_n = \alpha_{z(1)} \cdots \alpha_{z(n)}$. For any $q_1, q_1'$ lying
outside of $U_{z(1)}$, the cross-ratio
$(g_n q_n, g_n q_n' ; q_1, q_1')$ must be at least $n$. It follows
that (up to uniform additive and multiplicative error) the distance
from $g_n$ to the identity in $\Gamma$ is at least $n$.

\section{Proofs of Theorem \ref{thm:main} and Corollary~\ref{cor:cell}} \label{sec:proof} 

We can now prove the main theorem.  Recall our standing assumptions: $\Gamma$ is a hyperbolic group with fixed generating set $\mathcal{S}$.  We have fixed a metric $\dvis$ on $\bgam$, and the constant $D$ from Lemma~\ref{lem:distantpairs}.  Moreover in Definition~\ref{def:epsilon} 
we  fixed a neighborhood $\mathcal{U}$ of the identity in the space of continuous self-maps of $\bgam$, and a constant $\epsilon<D/4$ small enough so that $\mathcal{U}$ contains all maps $f$ such that $d(x, f(x)) \leq \epsilon$ for all $x \in \bgam$.

We have also fixed covers $\{W_z\supset V_z\}_{z\in I}$ in Definition~\ref{def:covers} that define a coding of boundary points as in Section \ref{sec:setup}, so that the results of Section \ref{sec:codingproperties} follow.  Recall also that the sets $W_z$ have diameter bounded by $\epsilon$; we will use this property later.

Our goal is to specify a neighborhood $\mathcal{V}$ of the standard boundary action in $\Hom(\Gamma, \Homeo(\bgam))$ such that every action in $\mathcal{V}$ is an extension of the standard boundary action via a semi-conjugacy in $\mathcal{U}$. 

\begin{definition} \label{def:same_comb}  We say that $\rho \in \Hom(\Gamma, \Homeo(\bgam))$ {\em has the same combinatorics as the standard boundary action} if the following hold for every $y,z\in I$.
\begin{enumerate}
\item $\overline{V}_z  \cap (\alpha_y^{-1}\overline{V}_y) \neq \emptyset$ iff $\overline{V}_z  \cap (\rho(\alpha_y)^{-1}\overline{V}_y) \neq \emptyset$, and 
\item if $\overline{V}_z  \cap (\alpha_y^{-1}\overline{V}_y) \neq \emptyset$, then $\overline{W_z} \subset \rho(\alpha_y)^{-1}(W_y)$. 
\end{enumerate}
\end{definition} 
Note that ``having the same combinatorics'' is an open
condition, because of Property~\ref{itm:cleanoverlap} from Definition~\ref{def:covers} of our covers.

\begin{definition}
  Suppose $\rho$ has the same combinatorics as the standard boundary action.
  We say that an infinite path $\{z(k)\}_{k\in\bN}$ in $\mc{G}$ is a \emph{$(\mc{G}, \rho)$--coding of $p$} if 
  \[ p \in \bigcap_{k = 0}^\infty\rho(\alpha_{z(1)})\cdots \rho(\alpha_{z(k)})\overline{W}_{z(k+1)}.\]
\end{definition}

Since Lemma \ref{lem:codingsexist} only used the intersection pattern of the sets $V_i$ and their images under the action of $\Gamma$, its proof applies verbatim to show the following.

\begin{lemma}  \label{lem:perturbedcodingsexist}
If $\rho \from \Gamma \to \Homeo(\bgam)$ has the same combinatorics as the standard boundary action, then every $p \in \bgam$ has a  $(\mc{G}, \rho)$--coding.
\end{lemma} 

Corollary~\ref{cor:unique} stated that for the standard action, each $\mc{G}$--coding determined a \emph{unique} point.  But the proof of that corollary used strongly that the action on $\bgam$ was induced by the isometric action on the Cayley graph.  Indeed there may be a nondegenerate closed subset of $\bgam$ all of whose points share the same $(\mc{G},\rho)$--coding.

\begin{definition}[The neighborhood $\mc{V}$]\label{def:V}
  We now describe a neighborhood of the standard boundary action in $\Hom(\Gamma, \Homeo(\bgam))$ that will satisfy the requirements of the theorem.  In the following description $F$ is the finite set from Lemma~\ref{lem:F} and  $N$ is the constant from the Uniform Contraction Lemma~\ref{lem:uniform}.  Similarly for $s\in S$, the set $F_s$ is the finite set from Lemma~\ref{lem:Fs}, and the constant $N_s$ is the constant from Corollary~\ref{cor:uniform_for_conj}.
  We let $\mc{V}$ be some neighborhood sufficiently small that all the following hold.
  \begin{enumerate}[label = (V\arabic*)]
  \item If $\rho\in \mc{V}$ then $\rho$ has the same combinatorics as the standard boundary action, in the sense of Definition~\ref{def:same_comb}.
  \item\label{itm:perturbnest} If $\rho\in \mc{V}$, $y,z\in I$, and $z(1)\ldots z(N)$ is a length-$N$ path in $\mc{G}$ so that $\alpha_{z(1)}\cdots \alpha_{z(N)} \overline{W}_z\subset f W_y$ for some $f\in F$, then
\begin{equation*}
  \rho(\alpha_{z(1)})\cdots \rho(\alpha_{z(N)})
  \overline{W}_z \subset \rho(f)W_{y} 
\end{equation*}
\item\label{itm:perturbnest_s} If $\rho\in \mc{V}$, $s\in \mc{S}$, $y,z\in I$, and $z(1)\ldots z(N_s)$ is a length-$N_s$ path in $\mc{G}$ so that $ \alpha_{z(1)} \ldots  \alpha_{z(N_s)}\overline{W}_{z}  \subset fs W_{y} $ for some $f\in F_s$, then
  \begin{equation*}
    \rho(\alpha_{z(1)}) \ldots  \rho(\alpha_{z(N_s)})\overline{W}_{z}  \subset \rho(fs) W_{y} .
  \end{equation*}
  \end{enumerate}
\end{definition}
\medskip
\putinbox{From here on we fix a representation $\rho$ from the neighborhood $\mc{V}$ just defined.}
\medskip
 
Our next goal is to define a function $\Phi\from \bgam \to 2^{\bgam}$ associating each point $p \in \bgam$ that is $\mc{G}$--coded
by a sequence $\{z(k)\}_{k\in\bN}$ to the closed set \[\bigcap_{k = 0}^\infty\rho(\alpha_{z(1)})\cdots \rho(\alpha_{z(k)})\overline{W}_{z(k+1)}.\]
The sets $\Phi(p)$ will be the fibers of our semi-conjugacy.
To show $\Phi$ is well defined, we 
use the following.

\begin{lemma} \label{lem:well-defined} 
If $\{z(k)\}_{k\in\bN}$ and $\{y(k)\}_{k\in\bN}$ are two distinct $\mc{G}$--codings of $p$, then 
\[ \bigcap_{k = 0}^\infty\rho(\alpha_{z(1)})\cdots \rho(\alpha_{z(k)})\overline{W}_{z(k+1)} = \bigcap_{k = 0}^\infty\rho(\alpha_{y(1)})\cdots \rho(\alpha_{y(k)})\overline{W}_{y(k+1)} \]
\end{lemma}
\begin{proof}
  It suffices to show that for any finite $k$, we can find some $n$
  such that
  \begin{equation}\label{eq:somenworks}
    \rho(\alpha_{z(1)})\cdots \rho(\alpha_{z(n)})\overline{W}_{z(n+1)}
    \subseteq \rho(\alpha_{y(1)})\cdots
    \rho(\alpha_{y(k)})\overline{W}_{y(k+1)}.
  \end{equation}
  Given $k$, we choose $n(k)$ as in Lemma~\ref{lem:F}, so that if
  $g_k = \alpha_{z(1)} \cdots \alpha_{z(k)}$ and
  $h_k = \alpha_{y(1)} \cdots \alpha_{y(k)}$, then $h_k = g_{n(k)}f$ for
  some $f \in F$.

  Then we choose $N$ as in Lemma
  \ref{lem:uniform}, so that
  \[
    \alpha_{z(1)}\cdots \alpha_{z(n(k) + N)}\overline{W}_{z(n(k) + N
      +1)} \subseteq \alpha_{y(1)}\cdots
    \alpha_{y(k)}W_{y(k+1)}.
  \]
  We claim that the containment above still holds even after we
  replace each $\alpha_k$ with its perturbed image $\rho(\alpha_k)$,
  i.e. that
  \[
    \rho(\alpha_{z(1)})\cdots \rho(\alpha_{z(n(k) +
      N)})\overline{W}_{z(n(k) + N +1)} \subseteq
    \rho(\alpha_{y(1)})\cdots
    \rho(\alpha_{y(k)})W_{y(k+1)}.
  \]
  In other words, \eqref{eq:somenworks} is satisfied with $n = n(k) + N$.
  To prove the claim, multiply each side of the above inclusion
  by $\rho(g_{n(k)})^{-1}$ to obtain one 
  of the conditions assumed in Item~\ref{itm:perturbnest} of Definition~\ref{def:V}.
\end{proof} 
Thus, the following gives a well-defined map from $\bgam$ to the space of closed subsets of $\bgam$.  

\begin{definition} \label{def:Phi}
  Let
  \[ \Phi(p) := \bigcap_{k = 0}^\infty\rho(\alpha_{z(1)})\cdots \rho(\alpha_{z(k)})\overline{W}_{z(k+1)},\] where $\{z(k)\}_{k\in\bN}$ is any coding of $p$.  
\end{definition} 

\begin{lemma}[Equivariance] \label{lem:equivariance} 
For any $p \in \bgam$ and $g \in \Gamma$, we have 
\[\Phi(gp) = \rho(g)(\Phi(p)).\]
\end{lemma}

This proof is where we use the fact that the nesting conditions for conjugates by generators from Item~\ref{itm:perturbnest_s} also hold under our perturbation.  

\begin{proof}[Proof of Lemma \ref{lem:equivariance}]
It suffices to prove the statement for an element in the finite generating set $\mathcal{S}$, then apply iteratively.  
Let $s \in \mathcal{S}$.   Let $\{z(k)\}_{k\in\bN}$ be any $\mc{G}$--coding of $sp$ and let $\{y(k)\}_{k\in\bN}$ be a $\mc{G}$--coding of $p$.  
The left-hand side of the claimed equality is, by definition, $\Phi(sp) = \bigcap_{k = 0}^\infty\rho(\alpha_{z(1)})\cdots \rho(\alpha_{z(k)})\overline{W}_{z(k+1)}$. 

The right-hand side is 
\begin{align*}
\rho(s)(\Phi(p)) &= \rho(s) \bigcap_{k = 0}^\infty\rho(\alpha_{y(1)})\cdots \rho(\alpha_{y(k)})\overline{W}_{y(k+1)}\\
&= \bigcap_{k = 0}^\infty\rho(s \alpha_{y(1)} s^{-1})\cdots \rho(s \alpha_{y(k)} s^{-1}) \rho(s) \overline{W}_{y(k+1)}.
\end{align*}
Since $\{y(k)\}_{k\in\bN}$ gives a $\mc{G}^s$--coding of $sp$ we may apply
Corollary \ref{cor:uniform_for_conj} to find some $N_s$ such that for
all $k$,
\[g_{n(k)+N_s} \overline{W}_{z(n(k)+N_s+1)} \subset h_{k}
  sW_{y(k+1)},\] where $g_k$, $h_k$, and $n(k)$ are as in Lemma~\ref{lem:Fs}.
Explicitly, this means that for any $k$,
\[ \alpha_{z(1)} \ldots \alpha_{z(n(k)+N_s)} \overline{W}_{z(n(k)+N_s+1)}  \subset s \alpha_{y(1)} s^{-1} \cdots s \alpha_{y(k)} s^{-1} sW_{y(k+1)}, \]
which gives, after multiplying both sides on the left by $(\alpha_{z(1)}\ldots \alpha_{z(n(k))})^{-1}$,  one of the finitely many containments
\[  \alpha_{z(n(k)+1)} \ldots  \alpha_{z(n(k)+N_s)}\overline{W}_{z}  \subset fs W_{y} \]
with $f\in F_s$ as in \eqref{eq:finite_s}.

Our assumption~\ref{itm:perturbnest_s} implies that this containment is
still satisfied after perturbation, i.e.
\[  \rho( \alpha_{z(n(k)+1)}) \ldots  \rho(\alpha_{z(n(k)+N_s)})\overline{W}_{z}  \subset \rho(fs) W_{y}. \]  After multiplying on the left by $\rho(\alpha_{z(1)}\ldots \alpha_{z(n(k))})$ we obtain
\[\rho(\alpha_{z(1)} \cdots \alpha_{z(n(k) + N_s)})
  \overline{W}_{z(n(k)+N_s+1)} \subset \rho(s\alpha_{y(1)}s^{-1}
  \cdots s\alpha_{y(k)}s^{-1}) \rho(s)W_{y(k+1)}\]
which shows that $\Phi(sp) \subset \rho(s)(\Phi(p))$. Applying the
same argument using $s^{-1}$ we also see that
\[
  \rho(s)(\Phi(p)) = \rho(s)\Phi(s^{-1}sp) \subset
  \rho(s)\rho(s^{-1})\Phi(sp) = \Phi(sp).
\]
\end{proof}

Combined with Lemma~\ref{lem:codingsexist} which states that every point has a $(\mathcal{G}, \rho)$ coding, the next lemma shows that the sets $\Phi(p)$ partition $\bgam$ as $p$ ranges over $\bgam$.
\begin{lemma} \label{lem:nonintersection}
  If $p \neq q$, then $\Phi(p) \cap \Phi(q) = \emptyset$.
\end{lemma}

\begin{proof}
First consider the case where $\dvis(p, q) > D$.  Then for any coding $\{z(k)\}_{k\in\bN}$ of $p$ and $\{y(k)\}_{k\in\bN}$ of $q$, respectively, we have $W_{z(1)} \cap W_{y(1)} = \emptyset$.  Since $\Phi(p) \subset W_{z(1)}$ and $\Phi(q) \subset W_{y(1)}$, this proves the lemma in this case.  
Lemma \ref{lem:distantpairs} (that any pair $a \neq b$ can be taken to a pair separated by distance $D$ by some group element) and \ref{lem:equivariance} (equivariance) reduce the general case to this one.  
\end{proof}

In summary, Lemma \ref{lem:perturbedcodingsexist}, Lemma \ref{lem:well-defined},
and Lemma \ref{lem:nonintersection} together imply that the sets
$\Phi(p)$ give a partition of $\bgam$, indexed by the points in
$\bgam$.  Thus this partition
defines a surjection
$\phi:\bgam \to \bgam$, determined by the condition
\[
  \phi(x) = p \iff x \in \Phi(p).
\]
Lemma~\ref{lem:equivariance} implies the function $\phi$ is \emph{$\rho$--equivariant} in the sense that for every $g\in \Gamma$ and $x\in\bgam$,
\begin{equation*}
  g\phi(x) = \phi(\rho(g)x).
\end{equation*}

\begin{lemma} \label{lem:phi_near_id}
  For every $x\in \bgam$, we have $\dvis(x,\phi(x))\le \epsilon$.
\end{lemma}
\begin{proof}
  Let $p = \phi(x)$, equivalently $x \in \Phi(p)$.  For any $\mc{G}$--coding $\{z(k)\}_{k\in\bN}$ of $p$,
  \[  \Phi(p) = \bigcap_{k = 0}^\infty\rho(\alpha_{z(1)})\cdots \rho(\alpha_{z(k)})\overline{W}_{z(k+1)}. \]
  In particular $p$ and $\Phi(p)$ are both contained in $\overline{W}_{z(1)}$, which has diameter at most $\epsilon$.  Since $x\in \Phi(p)$, $\dvis(x,p)\le \epsilon$.
\end{proof}

We have now verified all of the conditions needed for $\phi$ to be a
semi-conjugacy in the specified neighborhood $\mathcal{U}$ of
the identity, except for the fact that $\phi$ is continuous. This
last condition is implied by the properties already established, as follows.
\begin{lemma} \label{lem:continuous}
  Let $\rho:\Gamma \to \Homeo(\bgam)$ be an action of $\Gamma$ on its
  Gromov boundary.  For any $\epsilon_0 < D/4$, if $\phi:\bgam \to \bgam$
  is a $\rho$--equivariant surjection satisfying
  $\dvis(x, \phi(x)) < \epsilon_0$ for all $x \in \bgam$, then $\phi$ is
  continuous.
\end{lemma}
\begin{proof}
  We proceed by contradiction, and suppose that for a sequence
  $\{x_n\}_{n \in \bN}$ in $\bgam$, $x_n$ converges to $x$, but
  $\phi(x_n) = p_n$ does not converge to $\phi(x) = p$. Taking a
  subsequence, we may assume that $p_n$ converges to $q \ne p$. Since
  $\phi$ is $\rho$--equivariant, we may use Lemma
  \ref{lem:distantpairs} to assume that $\dvis(p, q) > D$.
  
  For sufficiently large $n$ we have $\dvis(x_n, x) < \epsilon_0$. Then
  by the triangle inequality we have
  $\dvis(p_n, p) \le \dvis(p_n, x_n) + \dvis(x_n, x) + \dvis(x, p) <
  3\epsilon_0$ and thus $\dvis(p_n, q) > D - 3\epsilon_0 > D/4$,
  contradicting the fact that $p_n \to q$.
\end{proof}

This concludes the proof of Theorem \ref{thm:main}.

\begin{proof}[Proof of Corollary~\ref{cor:cell}]
A cellular set (see \cite[II.6]{Daverman}) is a nested intersection of closed cells; a cellular \emph{map} is one whose fibers are cellular.   Corollary~\ref{cor:cell} states that in the case $\bgam$ is a topological sphere, the semi-conjugacy $h$ that we construct is cellular and can be approximated by homeomorphisms. 

By Remark~\ref{rem:cell}, the sets $\overline{W}_z$ can be taken to be closed cells.  For any $p\in S^n$, the fiber $h^{-1}(p)$ is equal to the set $\Phi(p)$ from Definition~\ref{def:Phi}.  Since our perturbed action has the same combinatorics as the standard action, this is a nested intersection of closed cells, with each closed cell contained in the interior of the previous one.  
Thus by definition the semi-conjugacy $h$ is cellular.

A cellular map is also \emph{cell-like} meaning that each of its fibers has the shape of a point.  Theorem \ref{thm:aqs} below implies that $h$ can be approximated by homeomorphisms, which concludes the proof of the corollary. 
\end{proof} 

\begin{theorem}\cite{Armentrout,Quinn,Siebenmann} \label{thm:aqs}
  A cell-like surjection of closed $n$--manifolds can be approximated by homeomorphisms.
\end{theorem}
\begin{proof}
  When $n=1$ the only possibilities for cell-like sets are intervals, and the result is easy to prove.  For $n=2$ see \cite[IV.25, Corollary IA]{Daverman}.  Armentrout  \cite{Armentrout} proved the theorem in case $n=3$ modulo the Poincar\'e Conjecture, which was proved by Perelman~\cite{MorganTian}.  The case $n=4$ is due to Quinn~\cite{Quinn}.  The higher dimensional cases are covered by a theorem of Siebenmann~\cite{Siebenmann} (see also \cite[IV.24, Corollary 3A]{Daverman}).
\end{proof}

\section{Global stability for groups with $S^1$ boundary}\label{sec:circle}

In the special case of groups with boundary $S^1$, namely, fundamental groups of closed hyperbolic orbifolds and surfaces, our techniques give a new proof of the main result of \cite{Mann} stated below, generalizing it to groups with torsion as well. Matsumoto \cite{MatsumotoBP} gives an alternate proof in the torsion-free case, which draws inspiration from Markov partitions and combination theorems.   The proof we sketch in this section is a different, perhaps simpler approach also inspired by symbolic dynamics.

We recall the statement of Theorem \ref{thm:global}:
\globalstab*

We fix a group $\Gamma$ and an action $\rho_0$ as in the statement.  Let $\rho_\partial\from \Gamma\to \Homeo(S^1)$ be the standard boundary action.  By assumption, there is a degree $k$ covering map $\psi\from S^1\to S^1$ so that $\rho_\partial(\gamma)\circ \psi = \psi \circ \rho_0(\gamma)$ for all $\gamma\in \Gamma$.

Let $C$ be the set of actions semi-conjugate to $\rho_0$ by a monotone, degree one semi-conjugacy.  Our goal is to show that $C$ is connected, closed, and open.  We will recall some results of Ghys~\cite{ghys,ghys_survey} and Matsumoto~\cite{Matsumoto86} to show that $C$ is connected and closed.  The techniques of the current paper will then be used to show that $C$ is open.

\begin{lemma}
  $C$ is closed.
\end{lemma}
\begin{proof}
  Thinking of $S^1$ as $\bR/\bZ$, any homeomorphism $f$ can be lifted (in various ways) to a homeomorphism $\tilde{f}\from \bR\to \bR$.  The \emph{translation number} \[\trot(\tilde{f}) = \lim_{n\to\infty}\frac{f(0)-0}{n}\] is only well-defined up to an integer, but for a pair $f,g\in \Homeo(S^1)$, the difference
  \[ \tau(f,g) = \trot(\tilde{f}\tilde{g}) - \trot(\tilde{f}) - \trot(\tilde{g}) \]
  is well-defined, and varies continuously in $f$ and $g$. 
 Following work of Ghys, Matsumoto showed that the semiconjugacy class of an action on $S^1$ is determined by the rotation numbers ($\trot(\tilde{f})$ reduced modulo $\bZ$) of the generators together with the values of $\tau$ on $\Gamma\times\Gamma$ \cite{Matsumoto} (see also \cite[5.11]{MannHandbook}).  In particular the condition of \emph{not} being semi-conjugate to a given action is clearly open.
\end{proof}

\begin{lemma}\label{lem:pathconnect}
  $C$ is path connected.
\end{lemma}
\begin{proof}
  Let $\rho\from \Gamma\to \Homeo(S^1)$ be semi-conjugate to $\rho_0$ by a degree one monotone map.  Ghys showed in \cite[Th\'eor\`eme A]{ghys} that this implies that $\rho$ has the same \emph{bounded Euler class} as $\rho_0$. (This is the theorem refined by Matsumoto's result referred to above -- the function $\tau$ is a representative of this class.)
  If the action of $\rho$ on $S^1$ is minimal, then
  \cite[Theorem 6.5]{ghys_survey} implies that $\rho$ is actually conjugate to $\rho_0$ (which is itself minimal) by an orientation preserving homeomorphism.
  Since $\Homeo_+(S^1)$ is path connected, such an action can be connected by to $\rho_0$ by a path of such conjugacies.

  If on the other hand $\rho$ is not minimal, then \cite[Proposition 5.6]{ghys_survey} implies that there is an invariant Cantor set $K\subset S^1$ in which every orbit is dense.  Collapsing the components of $\overline{S^1\minus K}$ gives a semi-conjugacy taking $\rho$ to an action $\overline{\rho}$ with a dense orbit, which can be connected to $\rho_0$ as above.  
  But this collapse can be realized as a limit of a path of conjugacies taking the components of $\overline{S^1\minus K}$ to shorter and shorter intervals in $S^1$.
\end{proof}

\begin{lemma} \label{lem:conjtoU}
Given any neighborhood $\mathcal{U}$ of $\rho_0$, any action semi-conjugate to $\rho_0$ by a degree one monotone map can be conjugated into $\mathcal{U}$. 
\end{lemma} 
\begin{proof}
  As in the last paragraph of the proof of Lemma~\ref{lem:pathconnect}, we suppose that $\rho$ is semi-conjugate to $\rho_0$ by a monotone degree one map, but not conjugate.  We saw there that there is a path $h_t$ of semiconjugacies for $t\in [0,1]$ so that $h_0$ is the identity, $h_1$ collapses $\rho$ to an action $\overline{\rho}$ which is conjugate to $\rho_0$, and $h_t$ is a homeomorphism for $t<1$.  Let $c_t$ be a path of conjugacies starting at the identity, and so that $c_1$ conjugates $\rho$ to $\rho_0$.  We see that the path $c_t\circ h_t$ is a path of semi-conjugacies continuously deforming $\rho$ to $\rho_0$.  These are conjugacies for $t<1$, establishing the lemma.
\end{proof}

We have shown that $C$ is closed and path connected.  To finish the proof of Theorem~\ref{thm:global} we need to show $C$ is open.  This is where the techniques of the present paper come into play.

\begin{proposition} 
$C$ is open.
\end{proposition} 

\begin{proof} 
 Theorem \ref{thm:main} states that there is an open neighborhood of the standard boundary action $\rho_\partial$ consisting of semiconjugate actions by cell-like maps $h$ that can be taken close to the identity.  Such maps of the circle are necessarily degree one and monotone.  

 In the case of a lift $\rho_0$ satisfying
 $\psi \rho_\partial(\gamma) = \rho_0(\gamma)\psi$, one may run the
 same strategy of proof to show $C^0$ stability of $\rho_0$ as
 follows. We use an automaton which is combinatorially the same as the
 one we construct in our original proof for $\rho_\partial$, but we
 replace the sets $V_z$ and $W_z$ with their preimages under
 $\psi$. The elements $\alpha_z \in \Gamma$ are the same as before.
 By Remark~\ref{rem:cell}, we may assume each set $W_z$ is the union
 of $k$ intervals whose closures are disjoint.  A coding given by our
 automaton now defines a $k$-tuple of points in $S^1$, and Lemma
 \ref{lem:quasigeodesic} now shows that a coding that contains a point
 $p$ gives a sequence in $\Gamma$ uniformly close to a geodesic ray
 from the identity to $\psi(p)$. The remaining results of Section
 \ref{sec:codingproperties} carry over.

We now say that an action $\rho$ {\em has the same combinatorics as $\rho_0$} if the conditions from Definition \ref{def:same_comb} hold not only for the sets $V_z$ and $W_z$, but for individual connected components.  Precisely, for any connected components $A_z, A_y$ of $V_z, V_y$ and $B_z, B_y$ of $W_z, W_y$ respectively, we require 
\begin{enumerate}
\item $\overline{A}_z  \cap (\rho_0(\alpha_y^{-1})\overline{A}_y) \neq \emptyset$ iff $\overline{A}_z  \cap (\rho(\alpha_y)^{-1}\overline{A}_y) \neq \emptyset$, and 
\item if $\overline{A}_z  \cap (\rho_0(\alpha_y^{-1})\overline{A}_y) \neq \emptyset$, then $\overline{B_z} \subset \rho(\alpha_y)^{-1}(B_y)$. 
\end{enumerate}

Choose a neighborhood $\mathcal{U}$ of $\rho_0$ small enough so that
elements of $\mathcal{U}$ have the same combinatorics as $\rho_0$, and
conditions \ref{itm:perturbnest} and \ref{itm:perturbnest_s} from
Definition \ref{def:V} hold.  The definition of $\Phi$ and the proof
of Lemmas \ref{lem:well-defined}, \ref{lem:equivariance},
\ref{lem:nonintersection} and \ref{lem:phi_near_id} carry over to give
an equivariant map $\phi$ that associates to each $k$-tuple of the
form $\psi^{-1}(x)$ a subset $\Phi(x) \subset S^1$ consisting of $k$
connected components, with the property that each connected component
of $\Phi(x)$ is Hausdorff distance at most $\epsilon$ from a point of
$\psi^{-1}(x)$.  Associating a connected component of $\Phi(x)$ to its
nearest point now gives the desired semi-conjugacy $\phi$ from $\rho$ to $\rho_0$.  

Lemma \ref{lem:continuous} can be adapted to this setting using that
$\rho_0$ is a lift of $\rho_\partial$ to a cover (and thus by
equivariance, we can act by $\rho_0$ on distinct $\phi$-images so that
they are distance at least $D$ apart in $\bgam$) and shows that this
semi-conjugacy is continuous.

We conclude that $\rho_0$ has a neighborhood $\mathcal{U}$ consisting of semi-conjugate actions, where the semi-conjugacy is given by a monotone, degree one map.  By Lemma \ref{lem:conjtoU}, an arbitrary semi-conjugate action $\rho$ can be {\em conjugated} by some $g \in \Homeo_+(S^1)$ to an action arbitrarily close to $\rho_0$, i.e. into $\mathcal{U}$.  Thus, $g^{-1}\mathcal{U}$ gives the desired neighborhood of $\rho$, proving that $C$ is open.  
\end{proof}

\newcommand{\etalchar}[1]{$^{#1}$}

\end{document}